\documentclass{article}
\usepackage{ijcai24}

\usepackage{times}
\usepackage{soul}
\usepackage{url}
\usepackage[hidelinks]{hyperref}
\usepackage[utf8]{inputenc}
\usepackage[small]{caption}
\usepackage{graphicx}
\usepackage{amsmath}
\usepackage{amsthm}
\usepackage{booktabs}
\usepackage{algorithm}
\usepackage{algorithmic}
\usepackage[switch]{lineno}
\usepackage{amssymb}
\usepackage{tabularx}
\usepackage{multirow}
\usepackage{threeparttable}

\usepackage{xcolor}

\title{An Adaptive Hybrid Genetic and Large Neighborhood Search Approach for Multi-Attribute Vehicle Routing Problems}

\author{
Weiting Liu$^{1}$
\and
Yunqi Luo$^{1,2}$\and
Yugang Yu$^{1}$\thanks{Corresponding author.}
\affiliations
$^1$School of Management, University of Science and Technology of China, China
$^2$International Institute of Finance, University of Science and Technology of China, China
\emails
weitinglau1999@mail.ustc.edu.cn,
yunqi\_luo@ustc.edu.cn,
ygyu@ustc.edu.cn
}
\begin{document}

\maketitle

\begin{abstract}

    % The adaptive large neighborhood search (ALNS) has shown considerable success in solving complex combinatorial optimization problems (COPs). ALNS select different heuristics adaptively to search neighborhood of a feasible solution for finding good solutions for optimization problems. However, for some complex COPs, ALNS may find a local optimal solution early. To address this limitation, we propose an adaptive hybrid genetic search and large neighborhood search (AHGSLNS) that selects variant individuals adaptively to avoid falling into local optimum early. The proposed method aims to help potential feasible solutions to find a good solution. We evaluate the proposed method on a multi-attribute vehicle routing problem, which is a classical COP, used in semi-automated storage and retrieval systems in the real world. The results show that our approach is competitive with vanilla ALNS. In addition, it obtained better performance in terms of convergence, convergence rate and stability. The implementation of our approach will be made publicly available.

   Known for its dynamic utilization of destroy and repair operators, the Adaptive Large Neighborhood Search (ALNS) seeks to unearth high-quality solutions and has thus gained widespread acceptance as a meta-heuristic tool for tackling complex Combinatorial Optimization Problems (COPs). However, challenges arise when applying uniform parameters and acceptance criteria to diverse instances of the same COP, resulting in inconsistent performance outcomes. To address this inherent limitation, we propose the Adaptive Hybrid Genetic Search and Large Neighborhood Search (AHGSLNS), a novel approach designed to adapt ALNS parameters and acceptance criteria to the specific nuances of distinct COP instances. Our evaluation focuses on the Multi-Attribute Vehicle Routing Problem, a classical COP prevalent in real-world semi-automated storage and retrieval robotics systems. Empirical findings showcase that AHGSLNS not only competes effectively with ALNS under varying parameters but also exhibits superior performance in terms of convergence and stability. In alignment with our dedication to research transparency, the implementation of the proposed approach will be made publicly available.
    
\end{abstract}

\section{Introduction}

    As a widely studied subject in operations research and computer science research communities, Combinatorial Optimization Problems (COPs) represent of a class of computational challenges that arise in various domain, encompassing logistics, scheduling and telecommunications, among others. The inherent complexity of COPs stems from their discrete nature and the exponential growth of potential solutions with the increasing size of problem. In addressing these NP-hard problems, the utilization of meta-heuristic algorithms stands as a testament to their efficacy in exploring the vest solution space and converge to near-optimal solutions.

    The Large Neighborhood Search (LNS) \cite{shaw1998using} is a well-regarded meta-heuristic that has proven effective in solving complex COPs.Grounded in the destroy-and-repair principle \cite{schrimpf2000record}, LNS systemically refines solutions by iteratively applying destroy and repair methods through heuristics. This approach facilitates the exploration of superior candidate solutions, following a promising search trajectory. Extending the original framework, the Adaptive Large Neighborhood Search (ALNS) \cite{ropke2006adaptive,wouda2023alns} allows for the concurrent utilization of multiple destroy and repair methods within the same search. 
    In ALNS, operators are dynamically assigned weights based on the corresponding performance while successful methods receive higher weights, enhancing their selection probability in subsequent iterations. This adaptive framework aims to deploy the most effective destroy and repair methods in each search iteration, potentially leveraging the strengths of multiple heuristics to optimize solutions for complex COPs. The inclusion of an acceptance criterion, such as simulated annealing, further enhances ALNS performance by considering a range of candidate solutions. The algorithm's flexibility is underscored by the influence of various parameters, including the severity of destroy methods, which shape distinct search trajectories.
    
    Despite the effectiveness of ALNS in solving COPs, determining optimal parameters and acceptance criteria for different instances of the same COP remains challenging. Variations in ALNS configurations, particularly in acceptance criteria, lead to diverse performances \cite{santini2018comparison}.
    To address these challenges, we propose an integration of ALNS with Genetic Search (GS) to obtain stable and high-quality solutions across different instances of the same variant of COPs by utilizing parallel mechanism which is able to allocate various acceptance criteria on ALNS\cite{liu2022data}. As a classical form of COPs, our focus lies on the Vehicle Routing Problem (VRP)\cite{dantzig1959truck}, more specifically, the Multi-Attribute Vehicle Routing Problem (MAVRP))\cite{vidal2014unified}. This variant of VRP encompassing multiple trips\cite{fleischmann1990vehicle,cattaruzza2016multi}, multiple depots\cite{renaud1996tabu,escobar2014hybrid,stodola2020hybrid}, backhauls\cite{kocc2018vehicle,cuervo2014iterated,queiroga2020exact}, 
    mixed backhauls\cite{liu2013heuristic}, 
    open route\cite{schrage1981formulation,sariklis2000heuristic,li2007open} and others, which is applied to the real-world task scheduling issues within the Carton Transfer Unit System (CTU system). 

    \begin{figure}[htbp!]
        \centering 
        \includegraphics[scale=0.33]{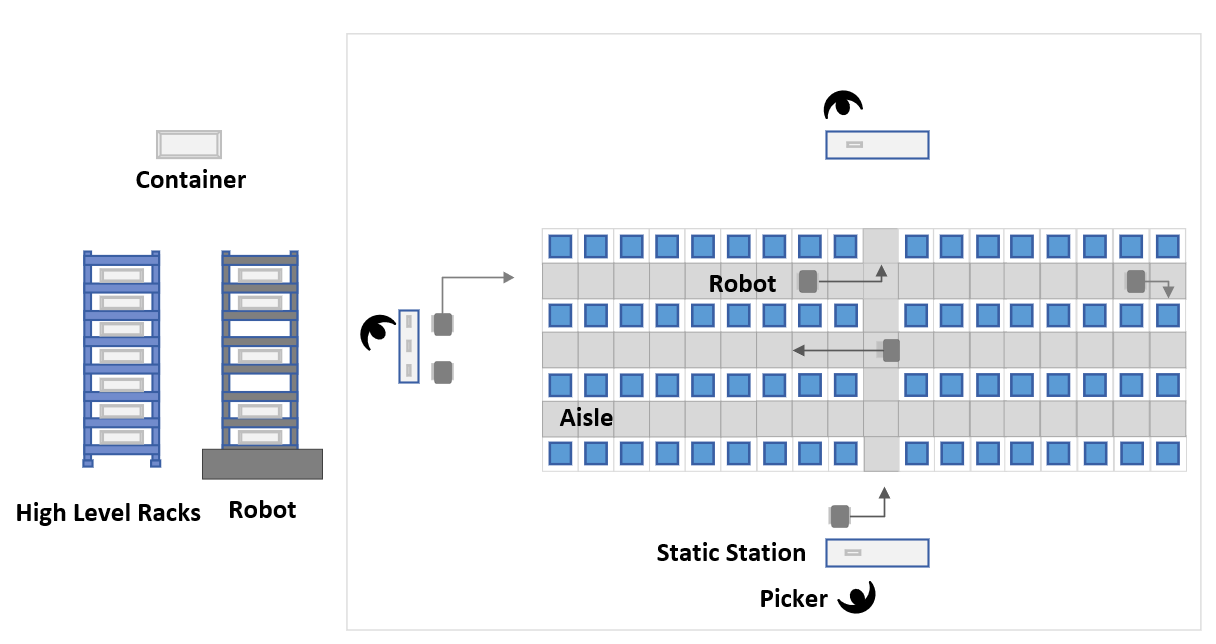} 
        \caption{System and schematic layout of a warehouse section with three static stations (See more details in https://www.hikrobotics.com/en/mobilerobot/CTU.)} 
        \label{fig:system}
    \end{figure} 
    
    The CTU system, as depicted in Figure 1, is a semi-automatic system designed for the storage and retrieval of production accessories and finished goods within the warehouse. The system comprises a central control system, containers, high-level racks, CTU robots equipped with multiple bins, and static stations. The central control system is responsible for task scheduling decisions, issuing instructions to CTU robots. Based on their assigned task types, such as storage tasks (linehaul customer tasks) or retrieval tasks (backhaul customer tasks), the robots autonomously transport and unload/load items between static stations (i.e., depots) and high-level racks along the aisles. Further, the SKUs sent to the static stations or returned to high-level racks undergo manual picking and storage processes conducted by workers. 
    
    Considering the automated loading and unloading features of the CTU robots in the system, in this paper, we extend the previous approach by integrating the search principles of ALNS and GS, constructing an adaptive meta-heuristic algorithm capable of adjusting parameters adaptively, named Adaptive Hybrid Genetic Search and Large Neighborhood Search (AHGSLNS). Our main contributions can be summarized as follows: (1) We formulate the mathematical model of the MAVRP problem in the CTU system under the automated loading and unloading features of CTU robots which also provide the benchmark for solving small-scale instances with the similar scenarios;(2) Utilizing an adaptive survival mechanism, we present a novel approach which selectively preserves elite individuals that are better suited for solving a specific instance when tackling different cases. This adaptive algorithm enables to tailor its performance to varying instances of same COP;(3) To further improve the performance of meta-heuristic algorithm on its convergence and stability, we use a cooperative evolution mechanism which is made up of the crossover and diversification operators to help the search to escape from the local optima trap at an earlier stage.

    We also conduct a series of experiments to verify the effectiveness of AHGSLNS. For small instances, the proposed algorithm can perfectly match the optimal results obtained from the exact model. And in large-scale counterpart, it surpasses ALNS by 100\% in terms of convergence and 73\% in stability. To further elucidate the performance of adaptive survival and cooperative evolution mechanisms, we conducted comparisons with AHGSLNS using various acceptance criteria across 15 instances. The experimental results indicate that, across 93.33\% of instances, AHGSLNS without cooperative evolution mechanism outperforms the average performance of 4 different ALNS variants while the AHGSLNS with cooperative evolution beats all the ALNS variants.

    \begin{figure}[htbp]
        \centering 
        \includegraphics[scale=0.35]{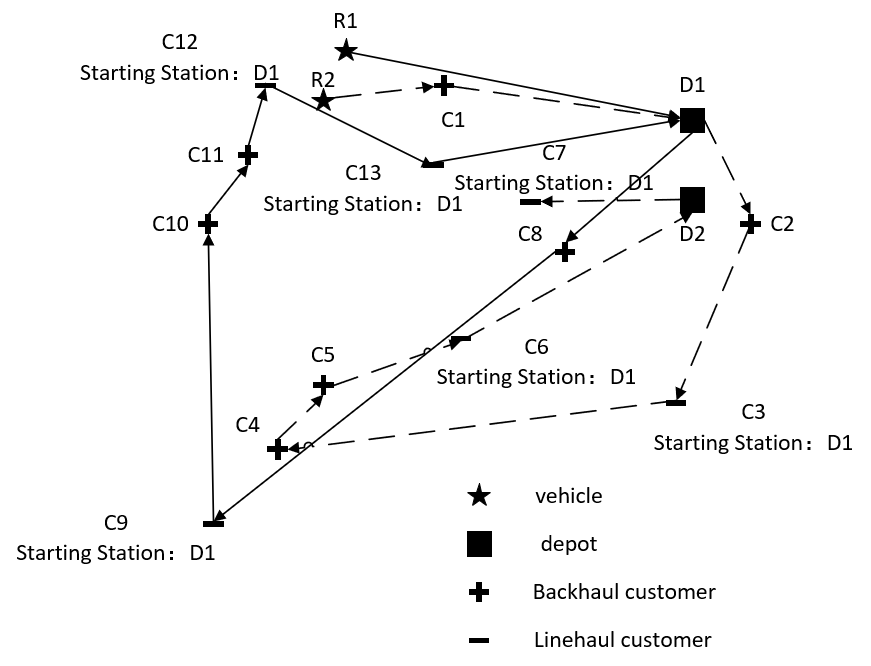} 
        \caption{Illustration of MAVRP in CTU system} 
        \label{fig:illustration}
    \end{figure}

\section{Problem Description and Mathematical Formulation}

\subsection{Problem Description}

    The MAVRP is formally defined on a directed graph $G=(N,A)$, where $N=D \cup K\cup L \cup B$ is the node set and $A=\{(i,j)|i,j \in N\}$ is the arc set. Here, $D=\{1,...,d\}$ embodies the depot (a.k.a. stationary station), whereas $K=\{d+m+n+1,...,d+m+n+a\}$ is a set of $a\in \mathbb{N_+}$ homogeneous vehicles with the capacity $Q\in \mathbb{N_+}$, $|K|=a$ and the vehicle node $k\in K$ represents the start terminals or the initial vehicle position.
    $L=\{d+1,...,d+m\}$ and $B=\{d+m+1,...,d+m+n\}$ represent $m$ linehaul customers and $n$ backhaul customers, respectively. Each node $i \in B \cup L$ is associated with a level of demand $q_i \in \mathbb{Z}\subseteq[-Q,Q]$. More specifically, the linehaul customers, hereafter referred to as linehauls, present a delivery demand denoted as $q_i>0,i\in L$, while backhaul customers, or backhauls, are characterized by a pickup demand, $q_i\leq 0, i \in B$.
    In addition, the MAVRP categorize routes into closed, open and return routes based on vehicle positions: (1) Closed routes refer to trajectories where both the start and end positions are located at the precise depot; (2) Open routes denote paths with start and end positions respectively situated at the depot and a designated customer location; and (3) Return trajectories characterize routes with start and end positions located at the vehicle's initial position and the depot, respectively. Any vehicle movement generated by one of the aforementioned routes is referred to as a ``trip" and we call a sequence of trips executed by the same vehicle as a ``journey"\cite{franccois2019adaptive}, denoted as $T=\{0,1,...,T_{max}\}$, where $T_{max}$ presents the maximum number of trips in each journey.We introduce the dummy node $u$ as the terminal node(i.e., the last customer) for each open trip. To enhance comprehension of MAVRP, Figure $2$ illustrates an instance with $d = 2$ depots, $k = 2$ vehicles, $m = 6$ linehauls and $n = 7$ backhauls.

\subsection{Integer Programming Formulation (ILP)}

    % The following assumptions are made: (1) Each customer, either linehaul or backhaul, is visited exactly once by exactly one trip of one vehicle; (2) To execute linehaul, a vehicle must depart from the depot node and then visit the linehaul node. (3) To execute backhaul, a vehicle must visit the backhaul node and finally arrive at the depot node; (4) Each open trip performed by the vehicles starts at the depot node and finishes at one of the customer nodes; (5) Each close trip performed by the vehicles starts at the vehicle nodes or depot node and finishes at the depot node; (6) A non-empty route must contain at least one linehaul customer or one backhaul customer; (7) The capacity $Q$ of the vehicle is not exceeded. We start by defining the following parameters and variables.

    Given the characteristics of the MAVRP,  we assume that (1) Each customer is visited exactly once by exactly one trip of one vehicle; (2) To execute linehaul, a vehicle must depart from the depot node and then visit the linehaul node; (3) To execute backhaul, a vehicle must visit the backhaul node and finally arrive at the depot node; (4) The capacity $Q=6$ of the vehicle is not exceeded. The goal of MAVRP is to schedule trajectories for each vehicle $k$ so that the maximum travel time of all vehicles can be minimized. Let $t_0$, $t_d$ respectively denote the index of the linehaul or backhaul customer's original or destination trips. We define $x_{ijk}^{(t_o)(t_d)}$ as the decision variable which equals to 1 if and only if a vehicle $k$ travel from customer $i$ to customer $j$ in trip $t_d$ where $t_o$ and $t_d$ are $j$'s original and terminal trip respectively. The loaded capacity of vehicle $k$ after visiting customer $i$ in the trip $t$ is represented as $U_{ik}^{(t)}\in\mathbb{Z}$ and the maximum travel time of all vehicles is donated by $C_{\max}$. We further introduce the nonnegative integer decision variable $e_i$ for customer $i$ for the related Miller-Tucker-Zemlin (MTZ) subtour elimination. The problem can be formulated as follows:

    \allowdisplaybreaks[4]
    \begin{align}
    %%%%%%%%%%%%%%%%%%%%%%%%%%%%%%%%%%%%%%%%%%%
    %%% (1) 最大完成时间目标函数
    %%%%%%%%%%%%%%%%%%%%%%%%%%%%%%%%%%%%%%%%%%%
    &\min_{x_{ijk}^{(t_o)(t_d)},U_{ik}^{(t)},e_i}C_{\max}\\
    \mbox{s.t.}\quad
    %%%%%%%%%%%%%%%%%%%%%%%%%%%%%%%%%%%%%%%%%%%
    %%% vertex constraint
    %%%%%%%%%%%%%%%%%%%%%%%%%%%%%%%%%%%%%%%%%%%
    %%% vehicle constraint
    %%% (2) 车辆入度约束
    %%% 第一次车辆节点不可以有入度
    &x_{ijk}^{(t_o)(t_d)} = 0, i \in N, t_o \in T, t_d \in T,j \in K, \notag \\
    &\hspace*{0.25in} k \in K \\
    %%% (3) 车辆出度约束
    %%% 车辆从车辆节点出发，可以到达depot或backhaul，但是只能出发一次
    &\sum_{j \in B \cup D} x_{kjk}^{(0)(0)} \le 1, k \in K \\
    %%% (4) 去程顾客入度约束
    %%% 车辆k只能从自己的节点出发，不能从其他的节点出发
    &x_{ijk}^{(0)(t_d)} = 0, i \in N, j \in L, t_d \in T, k \in K \\  
    %%% (5) 车辆去程顾客约束 
    %%% 车辆k只能从自己的节点出发，不能从其他的节点出发
    &\sum_{i \in K,i \ne k}\sum_{j \in N} x_{ijk}^{(0)(0)} = 0, k \in K \\ 
    %%% (6) 车辆出度位置限制约束
    %%% 第一次行程车辆不能前往去linehaul
    &x_{ijk}^{(0)(0)} = 0, i \in N, j \in L, k \in K \\ 
    %%% (7) 车辆出度行程限制约束
    %%% 除了第一次行程，车辆不能从车辆节点出发
    &x_{ijk}^{(t_o)(t_d)} = 0, i \in K,j \in N,t_o \in T \backslash \{0\} , \notag \\ 
    &\hspace*{0.25in} t_d \in T \backslash \{0\}, k \in K \\
    %%%%%%%%%%%%%%%%%%%%%%%%%%%%%%%%%%%%%%%%%%%
    %%% depot constraint
    %%% (8) 工作台第一次流量约束
    %%% 一辆车同一时间只能进入一个depot
    & x_{djk}^{(0)(0)} = 0,d \in D,j \in N, k \in K \\
    %%% (9) 工作台入度限制约束
    %%% 一辆车同一时间只能进入一个depot
    &\sum_{d \in D} \sum_{i \in N}x_{idk}^{(t_o)(t_d)} \le 1,k \in K,t_o \in T, t_d \in T \\
    %%% (10) 工作台非第一次行程出度限制约束
    %%% 除了第一次行程一个车辆同一时间最多从一个depot出发一次
    &\sum\limits_{d \in D} \sum_{j \in N}x_{djk}^{(t_o)(t_d)} \le 1,k \in K,t_o \in T \backslash \{0\}, \notag \\
    &\hspace*{0.25in} t_d \in T \backslash \{0\}  \\
    %%% (11) 工作台流量平衡约束
    %%% 上一刻进入depot的流量和下一刻从depot出发的流量需要保持平衡
    &\sum_{i \in N}x_{idk}^{(t_i)(t_i)} \ge \sum_{t_o \in \{1,...,t_i+1\}} \sum_{j \in B \cup L \cup D}x_{djk}^{(t_o)(t_i+1)}, \notag \\
    &\hspace*{0.25in} d \in D,k \in K, t_i \in T \backslash \{T_{max}\}  \\
    %%%%%%%%%%%%%%%%%%%%%%%%%%%%%%%%%%%%%%%%%%%
    %%% linehaul & backhaul constraint
    %%%%%%%%%%%%%%%%%%%%%%%%%%%%%%%%%%%%%%%%%%%
    %%% (12) 去程顾客和回程顾客入度限制约束
    %%% 去程顾客和回程顾客都需要有一个入度
    &\sum_{k \in K}\sum_{t_o \in T}\sum_{t_d \in \{t_o,...,T_{max}\}}\sum_{j \in N, i \ne j} x_{jik}^{(t_o)(t_d)} = 1, \notag  \\
    &\hspace*{0.25in} i \in B \cup L  \\
    %%% (13) 回程顾客第一次行程流平衡约束
    %%% 回程顾客哪次出发哪次到达,只有当次行程才可能为1
    %%% 回程顾客第一次行程需要保持流平衡，只有backhual和vehicle
    &\sum_{i \in B\cup K} x_{ijk}^{(0)(0)} = \sum_{i \in B\cup D} x_{jik}^{(0)(0)}, j \in B, k \in K \\
    %%% (14) 回程顾客非第一次行程流平衡约束
    %%% 回程顾客非第一次行程也需要保持流平衡
    &\sum_{i \in N}x_{ijk}^{(t_i)(t_i)} = \sum_{i \in N}x_{jik}^{(t_i)(t_i)},j \in B,t_i \in T, \notag  \\
    &\hspace*{0.25in} k \in K \\
    %%% (15) 去程顾客流平衡约束
    %%% 去程顾客需要和虚拟节点保持流平衡
    %%5 去程顾客必须有之前从某次行程工作台出发的流量，工作台->去程顾客，是终止行程
    &\sum_{t_1 \in \{0,...,t_d\}}\sum_{i \in B \cup D}x_{ijk}^{(t_1)(t_d)} = \sum_{t_2 \in \{0,...,t_d\}}  \sum_{i \in N \cup \{u\}} \notag \\
    &\hspace*{0.25in} x_{jik}^{(t_2)(t_d)}, j \in L,t_d \in \{1,...,T_{max}\}, k \in K \\
    %%% (16) 回程顾客解入度空间约束
    & x_{ijk}^{(t_o)(t_d)} = 0,  i \in N, j \in B, k \in K, t_o \in T, \notag \\
    &\hspace*{0.25in}  t_d \in T, t_o \ne t_d  \\
    %%% (17) 去程顾客解入度空间约束
    & x_{ijk}^{(t_o)(t_d)} = 0,  i \in N, j \in L, k \in K, t_o \in T, \notag \\
    &\hspace*{0.25in} t_d \in T, t_o \ne t_d  \\
    %%% (18) 顾客自环约束
    &x_{iik}^{(t_o)(t_d)} = 0, i \in B \cup L \cup K, t_o \in T, \notag \\
    &\hspace*{0.25in}  t_d \in T, k \in K \\
    %%%%%%%%%%%%%%%%%%%%%%%%%%%%%%%%%%%%%%%%%%%
    %%% 同质工作台场景约束
    %%%%%%%%%%%%%%%%%%%%%%%%%%%%%%%%%%%%%%%%%%%
    %%% 本次行程如果服务某个回程顾客必须返回任意一个工作台
    %%% (19) 回程顾客工作台限制约束
    %%% 回程顾客本次行程的顾客在本次行程最后的工作台都会被完成
    %%% 车辆k某次前往去程顾客，必须在这次前往某个工作台
    &x_{ijk}^{(t_i)(t_i)} = 1 \implies \sum\limits_{d \in D} \sum_{z \in B \cup L}x_{zdk}^{(t_i)(t_i)} = 1, \notag \\
    &\hspace*{0.25in} i \in  N\backslash K, j \in B, i \ne j, t_i \in T, k \in K \\
    %%% (20) 去程顾客出发约束
    &x_{ijk}^{(t_o)(t_d)} = 1 \implies \sum_{a \in N} \sum_{t_i \in \{1,...,t_d\}}x_{ad_jk}^{(t_o-1)(t_o-1)}=1, \notag \\
    &\hspace*{0.25in} i \in N, j \in L, i \ne j,   \notag \\
    &\hspace*{0.25in}  t_o \in T\backslash{0}, t_d \in \{t_o,...,T_{max}\}, k \in K \\
    %%% dummy end constraint
    %%% (21) 虚拟点限制约束
    %%% 车辆最多进入虚拟节点一次
    &\sum\limits_{i \in N}\sum_{t_o \in T}\sum_{t_d \in T}x_{iuk}^{(t_o)(t_d)} \le 1,k \in K\\
    %%% trip constraint
    %%% (22) 回程顾客和去程顾客装载量更新限制约束
    %%% 除了工作台更新其余点更新都要满足更新规则
    %%% 回程顾客位置更新，去程顾客终点位置更新
    &x_{ijk}^{(t_i)(t_i)} = 1 \implies U_{ik}^{(t_i)} + q_j = U_{jk}^{(t_i)}, i \in N,  \notag \\
    &\hspace*{0.25in} j \in B, i \ne j, t_i \in T, k \in K \\
    %%% (23) 回程顾客和去程顾客装载量更新限制约束
    %%% 回程顾客位置更新，去程顾客终点位置更新
    &x_{ijk}^{(t_o)(t_i)} = 1 \implies U_{ik}^{(t_i)} + q_j = U_{jk}^{(t_i)}, i \in N,  \notag \\
    &\hspace*{0.25in} j \in L, i \ne j, t_i \in T, t_o \in \{0,...,t_i\}, k \in K \\
    %%% (24) 车辆初始容量约束
    &U_{kk}^{(0)} = 0, k \in K \\
    %%% (25) 容量限制约束
    %%% 最大容量限制
    &U_{ik}^{(t)} \le Q, i \in N, t \in T, k \in K \\
    %%% 回程顾客在行程中完成，不需要考虑。去程顾客需要减去上次行程完成的，加上本次行程新加入的
    %%% (26) 工作台非第一次行程装载量更新约束
    %%% 回程顾客在行程中不会影响工作台的容量更新，因为路上接到，终点放下
    %%% 工作台装载量只需要更新去程顾客
    &x_{d_2jk}^{(t_1)(t_i)} = 1 \implies U_{d_2k}^{(t_i)} = \sum_{d_1 \in D \cup \{k\}} U_{d_1k}^{(t_i-1)}    \notag \\
    &+ \sum_{i \in N} \sum_{z \in L, d_z = d_2} \sum_{t_2 \in \{t_i,...,T_{max}\}} \sum_{j \in N}x_{izk}^{(t_i)(t_2)} \notag \\
    &+  - \sum_{t_1 \in \{0,...,t_i\}} \sum\limits_{z \in L} \sum_{i \in N} x_{izk}^{(t_1)(t_i-1)}, j \in N, d_2 \in D, \notag \\
    &\hspace*{0.25in}   t_i \in \{1,...,T_{max}\}, t_1 \in {0,...t_i}, k \in K \\
    %%% (27) 子环消除约束
    &e_{i} - e_{j} + 2 * Q * \sum_{k\in K}\sum_{t_o\in T}\sum_{t_d \in T}x_{ijk}^{(t_o)(t_d)} \le 2 * Q  \notag \\
    &\hspace*{0.25in} - 1, i\in B \cup L, j \in B \cup L, i \ne j \\
    %%%%%%%%%%%%%%%%%%%%%%%%%%%%%%%%%%%%%%%%%%%
    %%% makespan constraint
    %%%%%%%%%%%%%%%%%%%%%%%%%%%%%%%%%%%%%%%%%%%
    %%% (28) 最大完成时间约束
    &C_{\max} \ge \sum_{t_o \in T} \sum_{t_d \in T} \sum_{i \in N} \sum_{j \in N} x_{ijk}^{(t_o)(t_d)}, k \in K \\
    %%%%%%%%%%%%%%%%%%%%%%%%%%%%%%%%%%%%%%%%%%%
    %%% variable constraint
    %%%%%%%%%%%%%%%%%%%%%%%%%%%%%%%%%%%%%%%%%%%
    % decision variables
    %%% (29) 路径决策变量约束
    &x_{ijk}^{(t_o)(t_d)} \in \{0,1\}, i,j \in N, t_o \in T, t_d \in T, \notag \\
    &\hspace*{0.25in} k \in K \\
    %%% (30) 装载量决策变量约束
    &U_{ik}^{(t)} \ge 0, i \in N, t \in T, k \in K\\
    %%% (31) 子环消除非负变量约束
    & e_i \ge 0, i \in B \cup L \
    \end{align}
    
    Objective (1) aims to minimize the makespan(i.e., maximum travel time of all vehicles) of all vehicles. Constraints (2)-(7) are related to vehicle nodes. Constraint(2) restricts vehicles from entering the vehicle node during the first trip. Constraints (3)-(7) ensure that each vehicle can depart from the vehicle node and arrive at the backhauls or depot. 
    Constraints (8)-(11) pertain to depot node. Constraint (8)  ensures that vehicles do not depart from the depot during the first trip. Constraint (9) guarantees that each trip of vehicle which departs from other nodes (including the depot node) can enter depot only once. Constraint (10) represents that, apart from the initial trip, vehicles can leave the depot node at most once during each trip, then reach to any types of nodes. The depot-to-depot arcs with zero distance cost is considered to add the option of executing fewer than $T_{max}$ trips, these artificial trips will be called ``empty trips". Constraint (11) indicates that if a vehicle executes an open route then it cannot leave the depot node in the subsequent trips.
    (12)-(20) encapsulate constraints related to customer nodes.Constraint (12) ensures that each linehaul and backhaul must have an in-degree in order to be served. Constraints (13) and (14) are balance constraints for the first trip and remaining trips of backhaul nodes respectively, whereas constraint (15) enforces flow balance at linehaul customer nodes. Constraints (16) and (17) are variable feasibility constraints. Constraint (18)enforces a self-loop elimination constraint for customers. Constraint (19) is a back constraint, requiring vehicles to return to the depot node if a backhaul customer is served in a trip. Constraint (20) is a linehaul starting constraint, indicating that, if there is a linehaul node in a trip, the vehicle needs to start from the related station of the linehaul in previous trips. 
    Constraint (21) is related to a dummy node, ensuring that each vehicle has at most one open route.
    Constraints (22)-(26) are related to vehicle capacity. Constraints (22) and (23) ensures the updating load capacity of a vehicle if it travels in the trip. Constraint (24) is the initial load capacity constraint and (25) guarantees that the vehicle's capacity constraint is not violated. Constraint (26) is the updating constraint related to the capacity of depot node. Constraint (27) is the sub-tour elimination constraint. Equation (28) ensures that the maximum completion time of the system is greater than the maximum completion time of all vehicles. Constraints (29)-(31) bound the variables. It is worth noting that Equations (19), (20), (22), (23), and (26) are can be easily transformed  with Big-M method or directly solved by using Gurobi\cite{gurobi}.

\section{ALGORITHMS} 

    The framework of the AHGSLNS is described in Algorithm 1. Notice that the initialization phase in Lines 1-7 uses random constructive algorithm to generate initial population. In the adaptive survival phase, AHGSLNS uses PALNS (shown in Algorithm 2) to intensify the solution quality of individuals in population and eliminates worse individuals iteratively to adapt different instances (Lines 8-14). We use the parallel search mechanism to offset time consumption of population mechanism. Lines 15-26 describe the cooperative evolution phase which employs crossover and diversification operators to further improve the algorithm's performance.

    \begin{algorithm}[tb]
        \caption{Adaptive Hybrid Genetic Search and Large Neighborhood Search (AHGSLNS)}
        \label{alg:algorithm1}
        \textbf{Input}: $D,L,B,K$  //Depots set, Linehaul customers set, Backhaul customers set, Vehicles set.\\
        \textbf{Parameter}: $N, gen^{max},M,iter^{max}$ //Numbers of population individuals, Numbers of generation, Numbers of elite individuals, Max iterations of PALNS. \\
        \textbf{Output}: $x^b$ $//$The optimal solution. 
        
        \begin{algorithmic}[1] %[1] enables line numbers
            \STATE // Initialization. See section 3.2 for details.
            \STATE{P = \{\}}
            \WHILE{$|P| < N$}
                \STATE $s \gets$ random constructive algorithm  
                \STATE $\rho^{-}_{s}=(1,...,1);\rho^{+}_{s}=(1,...,1)$
                \STATE $x^{b}_{s}=x_{s}^c=x_{s}$,$P \gets P$ $\cup$ $\{s\}$
            \ENDWHILE
            % \STATE // Competitive elimination phase. 
            \STATE // Adaptive Survival. See section 3.3 for details.
            \WHILE{$|P| > M$}
                \STATE \textbf{parallel for} $s \in P$ \textbf{do}
                    \STATE$x^{b}_{s}$, $x_{s}^c$, $\rho^{-}_{s}$, $\rho^{+}_{s}$,$accept_{s}$,$DR_{s}^{set}$ $\gets$ PALNS($x^{b}_{s}$, $x_{s}^c$,  $\rho^{-}_{s}$, $\rho^{+}_{s}$,$accept_{s}$,$DR_{s}^{set}$,$iter^{max})$ 
                \STATE \textbf{end parallel for}
                \STATE{eliminate the worse individual $s_i, s_i \in P$}
            \ENDWHILE
            \STATE // Cooperative Evolution. See section 3.5 for details.
            \WHILE{$gens < gen^{max}$}
                \IF{crossover criterion is met}
                    \STATE Generate offspring from population. 
                \ENDIF
                \IF{diversification criterion is met}
                    \STATE Diversify population.
                \ENDIF
                \FOR{$s_i \in P$}
                    \STATE$x^{b}_{s}$, $x_{s}^c$, $\rho^{-}_{s}$, $\rho^{+}_{s}$,$accept_{s}$,$DR_{s}^{set}$ $\gets$ PALNS($x^{b}_{s}$,$x_{s}^c$,  $\rho^{-}_{s}$, $\rho^{+}_{s}$,$accept_{s}$,$DR_{s}^{set}$,$iter^{max})$ 
                \ENDFOR
            \ENDWHILE
            \STATE $x^b \gets$ find the best $x_{s}^b$ from $P$
            \STATE \textbf{return} solution $x^b$
        \end{algorithmic}
    \end{algorithm}
    
\subsection{Solution Representation, Individual Representation and Evaluation of Individual}

    A MAVRP solution $x$ is also represented by $a$ journeys which is consistent with the numbers of the vehicles. Each journey starts from a vehicle position and stop at the depot or one of the customers. We can transform feasible solution illustration in Figure 2 into journey solution representation as shown in Figure 3 where the characters within the parentheses indicate that the actual location of this customer is at the static station.

    \begin{figure}[ht!]
        \centering 
        \includegraphics[scale=0.2]{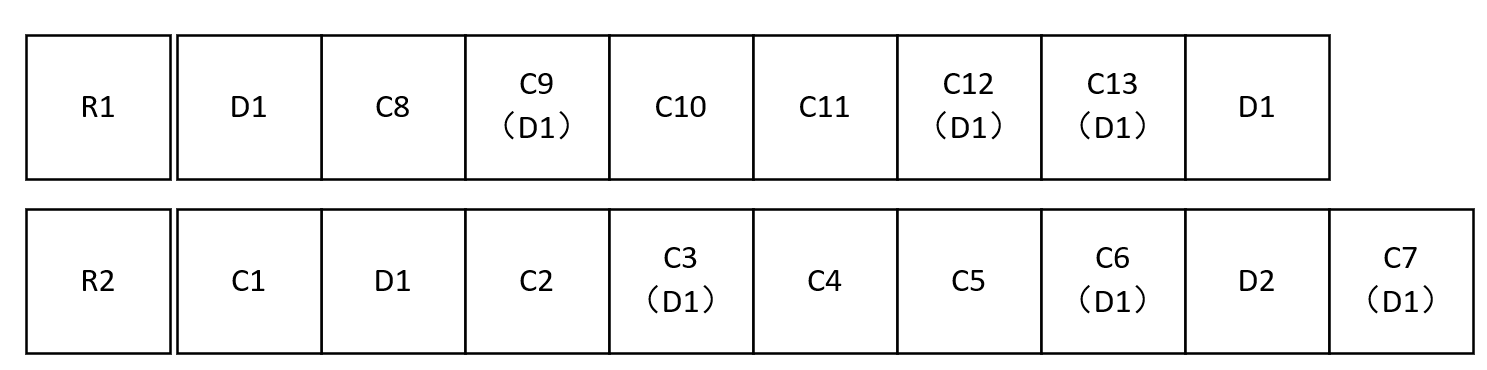} 
        \caption{Solution Representation} 
        \label{fig:representation}
    \end{figure} 

    Further, to evaluate the individual, we define function $f_k$ for all $k \in K$ which denotes the travel time of journey of the vehicle $k$ and let $MSP(x)$ be the expected cost makespan of a MAVRP solution $x$. The the fitness function of each solution is represented as
    \begin{align}
        fitness(x) = MSP(x) = \max\limits_{k \in K}\{f_k\}.
    \end{align}
    \noindent We define the individual $s$ as $(x_s^b$, $x_s^c$, $\rho^-$, $\rho^+$, $accept_{s}$, $DR_s^{set})$. In detail,  $x_s^b $ and $x_s^c$ are the best solution and current solution, respectively. The weights of destroy and repair methods are represented by $\rho^-$ and $\rho^+$ and the destroy-and-repair parameters set is given as $DR_s^{set}$ . We use $accept_{s}$ to denote the corresponding acceptance criteria. Note that fitness of the individual $s$ is the fitness of best solution $x_s^b$.

\subsection{Initialization Mechanism}

    To initialize a diversified population, AHGSALNS iteratively employs the random constructive algorithm multiple times to generate the individuals .

\subsubsection{Random Construction Algorithm}

    Note that the starting and ending locations of linehaul customers are predetermined, while the destination of backhaul customers requires decision-making. To create a feasible initial solution, we construct a three-phase method, which includes a customer assignment phase, a routing phase and a route assignment phase. Initially, linehaul customers are allocated to depot point and the closest backhaul customers are assigned to the depots that has the fewest customers. Subsequently, the savings algorithm\cite{clarke1964scheduling} is applied to generate closed trips for each depot, with the first customer of these trips selected randomly. Thirdly, closed routes are assigned to vehicles, and the vehicle journeys are constructed using a greedy policy to establish an initial feasible solution. Lastly, random ALNS parameters and an acceptance criterion are assigned to the corresponding individual.
    
    \begin{table*}[ht]
        \centering
        % \vspace{-2.0em}
        \label{tab:multi_simulation}
        \begin{tabular}{lrr r rrr r rrr}\toprule
        \multirow{ 2}{*}{Instances} &\multirow{ 2}{*}{BKS}&\multirow{ 2}{*}{ILP}& &\multicolumn{3}{c}{ALNS}&&\multicolumn{3}{c}{AHGSLNS}\\
        \cmidrule{5-7}\cmidrule{9-11}
          & &  &  &Best/Avg  &Gap& Std/CoV  && Best/Avg  & Gap & Std/CoV\\
        \midrule
        R\_2\_2\_5\_7 &\textbf{75.42} & \textbf{75.42} & &\textbf{75.42}/81.89 & 8.58\% & 4.78/5.84\% & &\textbf{75.42}/77.75 & 3.09\% & 1.38/1.77\% \\
        R\_2\_2\_8\_3 & \textbf{74.23} & \textbf{74.23} & &\textbf{74.23}/77.59 & 4.52\% & 4.26/5.49\% & &\textbf{74.23}/77.55 & 4.47\% & 4.08/5.26\% \\
        R\_2\_2\_7\_9& \textbf{84.13} & \textbf{84.13} & &\textbf{84.13}/87.95 & 4.54\% & 6.74/7.67\% & &\textbf{84.13}/84.88 & 0.89\% & 1.48/1.74\% \\
        R\_2\_2\_2\_10& \textbf{69.36} & \textbf{69.36} & &\textbf{69.36}/71.47 & 3.04\% & 1.05/1.48\% & &\textbf{69.36}/69.88 & 0.74\% & 1.04/1.49\% \\
        R\_2\_2\_12\_4 & - & 81.28 & &87.13/95.29 & 17.23\% & 5.19/5.45\% && 82.02/90.23 & 11.01\% & 3.33/3.69\% \\
        R\_2\_2\_3\_15 & -& 102.46 & &102.33/106.39 & 3.83\% & 2.65/2.49\% && 96.90/102.23 & -0.22\% & 3.48/3.41\% \\
        R\_2\_3\_9\_10 & -& 70.20 & &70.20/75.72 & 7.86\% & 3.62/4.78\% && 70.20/72.89 & 3.83\% & 1.51/2.08\% \\
        R\_2\_3\_12\_7 & -& 67.41 & &68.66/72.43 & 7.44\% & 2.39/3.30\% && 67.41/71.56 & 6.16\% & 2.67/3.77\% \\
        R\_2\_3\_17\_19 & -& 130.96 & &96.30/105.56 & -19.39\% & 6.43/6.09\% && 97.75/103.46 & -21.00\% & 3.62/3.51\% \\
        R\_2\_3\_21\_9 & -& 126.71 & &96.22/103.06 & -18.66\% & 5.90/5.72\% && 92.78/98.15 & -22.53\% & 3.86/3.93\% \\
        R\_2\_3\_19\_7\ & -& 101.66 & &88.39/93.19 & -8.33\% & 3.74/4.01\% && 85.84/90.34 & -11.13\% & 3.48/3.85\% \\
        R\_2\_3\_9\_21 & -& 130.85 & &92.02/99.59 & -23.88\% & 5.54/5.56\% && 88.28/95.87 & -26.73\% & 4.75/4.96\% \\
        R\_2\_3\_11\_16 & -& 93.57 & &77.52/88.39 & -5.53\% & 7.01/7.94\% && 83.87/86.57 & -7.48\% & 1.80/2.08\% \\
        R\_2\_3\_13\_13 & -& 90.84 & &81.35/85.51 & -5.87\% & 3.46/4.04\% && 77.26/80.91 & -10.93\% & 3.70/4.58\% \\
        R\_2\_3\_15\_16 & -& 141.55 & &85.69/94.16 & -24.08\% & 4.94/5.25\% && 84.70/89.94 & -36.46\% & 2.82/3.13\% \\
        \bottomrule
        \end{tabular}
        \caption*{Note: (a)``Best" and ``Avg" indicates the best and the average result of the algorithm when solving 10 times on a given instance, respectively; (b) Bold fold indicates that the current solution is optimal.}
        \caption{Comparison Results for Synthetic Instances among ILP, ALNS and AHGSLNS}
    \end{table*}

    \begin{table}
        \centering
        \begin{tabular}{crlcc}
            \toprule
            Type  & Set & Acceptance Criterion & ASM & CEM \\
            \midrule
            A1 & 1 & record-to-record travel & - & - \\
            A2 & 2 & record-to-record travel & - & - \\
            A3 & 3 & annealing acceptance & - & - \\
            A4 & 1 & hill climbing & - & - \\
            B1 & 1 & adaptive mechanism & \checkmark & -  \\
            B2 & 1 & adaptive mechanism & \checkmark & \checkmark \\
            \bottomrule
        \end{tabular}
        \label{tab:alns_parameters}
        \caption{Algorithm Setting for Variant of ALNS and AHGSLNS}
    \end{table}

    \begin{table*}[ht]
        \centering
        % \vspace{-2.0em}
        \label{tab:adaptive_maechism}
        \begin{tabular}{lrrrrrr}
            \toprule  
            Instance & Best/Avg(A1) &  Best/Avg(A2)  & Best/Avg(A3) & Best/Avg(A4)  & Best/Avg(B1)  & Best/Avg(B2)\\
            \midrule 
            R\_2\_2\_11\_3  & \underline{61.47}/67.13 & 68.29/69.68 & \underline{61.47}/65.32+ & 68.29/71.87  & \textbf{61.47}/\textbf{65.28} & \textbf{61.47}*/\textbf{63.40}* \\
            R\_2\_2\_7\_15  & 79.47/81.24+ & 79.17/84.49 & \underline{78.77}/81.49 & 79.55/82.80  & 79.17/\textbf{80.70} & 79.11*/\textbf{80.33}* \\
            R\_2\_2\_11\_11 & 76.54/81.87+ & \underline{76.45}/84.98 & 79.98/84.00 & 80.65/86.28 & 77.00/\textbf{79.91} & 77.00*/\textbf{79.71}* \\
            R\_2\_2\_15\_7 & 79.27/84.11+ & 79.92/84.88 & 80.46/84.70 & \underline{79.05}/87.13 & 79.55/\textbf{83.40} & \textbf{76.41}*/\textbf{79.72}* \\
            R\_2\_2\_16\_6 & 79.24/86.40 & 81.99/86.54 & \underline{78.86}/83.89+ & 83.77/86.79  & \textbf{78.32}/84.55 & \textbf{75.78}*/\textbf{81.34}* \\
            R\_2\_2\_12\_9 & \underline{74.29}/80.17+ & 78.21/82.76 & 79.80/83.52 & 79.89/84.11 & 78.40/\textbf{79.62} & 75.86*/\textbf{78.43}* \\
            R\_2\_2\_10\_11 & \underline{75.99}/80.04 & 77.18/78.43 & 77.85/81.19 & 76.22/77.81+ & 76.42/\textbf{77.37} & \textbf{75.35}*/\textbf{76.90}* \\
            R\_2\_2\_12\_11 & 89.71/91.85 & \underline{82.69}/90.46 & 87.09/90.37+ & 85.57/91.95 & 83.37/\textbf{88.60} & 85.17/\textbf{88.69} \\
            R\_2\_2\_17\_6 & 85.37/88.83+ & \underline{80.63}/90.94 & 83.95/89.30 & 87.65/94.22 & 81.28/\textbf{86.00} & \textbf{79.63}*/\textbf{88.04} \\
            R\_2\_2\_11\_13 & \underline{80.87}/84.53+ & 82.11/85.73 & 81.56/87.15 & 84.88/86.57 & \textbf{78.34}/\textbf{82.78} & \textbf{80.25}/\textbf{81.68}* \\
            R\_2\_2\_12\_10 & 79.70/84.50 & 82.60/84.93 & \underline{77.01}/83.05+ & 80.08/84.75 & 79.04/\textbf{82.55} &  \textbf{76.66}*/\textbf{80.62}* \\
            R\_2\_2\_6\_13 & \underline{79.54}/81.44+ & \underline{79.54}/82.36 & 80.85/84.33 & \underline{79.54}/81.75 & \textbf{79.54}/\textbf{80.64} & \textbf{79.54}*/\textbf{80.84} \\
            R\_2\_2\_9\_13 & \underline{76.40}/81.49+ & 78.54/82.49 & 77.29/81.90 & 79.55/82.44 & 77.54/\textbf{80.07} & 76.80*/\textbf{79.55}* \\
            R\_2\_2\_15\_18 & 105.35/110.32 & 109.83/115.05 & 103.36/115.50 & \underline{101.67}/108.52+ & 104.72/\textbf{106.74} & 101.69*/\textbf{105.40}* \\
            R\_2\_2\_19\_17 & 112.69/118.83+ & 110.56/120.91 & 119.03/126.09 & \underline{110.38}/125.33 & 110.63/\textbf{117.80} & \textbf{112.84}*/\textbf{118.09}* \\
            \bottomrule 
        \end{tabular}
        \caption*{Note:(a) Similarly to Table 1; (b) Underline and ``+" represent the best and average results among A1, A2, A3, and A4, respectively. (c) Bold font represents better results of B1 or B2 compared with that of ALNS's variants. (d) ``*" represents the better results among B1 and B2. }
        \caption{Performance of ALNS's Variants and AHGSLNS}
    \end{table*}

\subsection{Parameterized Adaptive Large Neighborhood Search (PALNS) Framework}

    % The procedure of PALNS is shown in pseudocode in Algorithm 2. 
    % Line 1 initialize the current iterations of PALNS. Lines 3-4 represents the process of adaptive neighborhood search using destroy methods and repair methods. Lines 2 and Lines 5-7 are stop criterion and acceptance criterion, respectively. Lines 8-10 are updating process of the best solution $x_s^b$. Line 11 is the adjustment process of $iter^c$ and weights of destroy methods and repair methods.

    % $\rho^{-}=(1,...,1);\rho^{+}=(1,...,1)$
    \begin{algorithm}[tb]
        \caption{Parameterized Adaptive Large Neighborhood Search (PALNS)}
        \label{alg:algorithm2}
        \textbf{Input}: $x_{s}^b, x_{s}^c,\rho_{s}^{-},\rho_{s}^{+},accept_{s},DR_{s}^{set},iter^{max}$  //Best solution of $s$, Current solution of $s$, Destroy method weights of $s$, Repair method weights of $s$, Acceptance criterion of $s$, Maximum number of iterations  \\
        \textbf{Output}: $x_{s}^b, x_{s}^c,\rho_{s}^{-},\rho_{s}^{+},accept_{s},DR_{s}^{set}$
        
        \begin{algorithmic}[1] %[1] enables line numbers
            \STATE $iter^{c}=0$   // current iteration
            \WHILE{$iter^{c} < iter^{max}$}
                \STATE select destroy and repair methods d $\in$ $\Omega^{-}$ and r $\in$ $\Omega^{+}$ using $\rho_{s}^{-}$ and $\rho_{s}^{+}$;
                \STATE $x_{s}^t=r(d(x_{s}^c))$
                \IF {$accept_{s}$($x_{s}^t, x_{s}^c$)}
                \STATE $x_{s}^c = x_{s}^t$
                \ENDIF
                \IF {c($x_{s}^t$) $<$ c($x_{s}^b$)}
                \STATE $x_{s}^b$=$x_{s}^t$
                \ENDIF
                \STATE update $iter^c$,$\rho_{s}^{-}$ and $\rho_{s}^{+}$
            \ENDWHILE
            \STATE \textbf{return} $x_{s}^b, x_{s}^c,\rho_{s}^{-},\rho_{s}^{+},accept_{s},DR_{s}^{set}$  
        \end{algorithmic}
    \end{algorithm}

\subsubsection{Destroy Methods}

    In view of the characteristics of the problem, we design a series of destroy operators, including the worst destroy method\cite{ropke2006unified}, related destroy method\cite{shaw1997new}, history destroy method\cite{ropke2006unified}, string destroy method\cite{christiaens2020slack}, and other problem-specific destroy methods. Each destroy method is equipped with parameters that regulate the severity and characteristics of solution destruction in the neighborhood search process, exemplified by the destroy severity parameter in the worst destroy method.
    
\subsubsection{Repair Methods}

    Similarly, a suite of repair operators are introduced, which encompass the basic greedy repair method\cite{ropke2006adaptive}, regret repair method\cite{potvin1993parallel}, random repair method\cite{ropke2006adaptive}, and other problem-specific repair methods. Repair methods also have parameters to control repairing characteristic of the solution, such as regret degree of regret repair method.

\subsubsection{Choosing a Destroy and Repair Method}

    Different destroy methods are applicable to distinct repair methods. To determine the heuristic selection, we assign weights to each heuristics with roulette wheel selection principle. Let $k_1$ represent the number of destroy methods with weights $w_i$ for $i \in \Omega^{-}$, $ \Omega^{-} = \{1,2,...,k_1\}$. The destroy method $d$ is selected with probability $\frac{w_{d}}{\sum_{z=1}^{k_1}w_z}$. The choice of repair methods bears similarity to that of destroy methods.

    %Similarly, $k_2$ repair methods with weight $w_j, j \in \Omega^{+}, \Omega^{+} = \{k_1+1,k_1+2,...,k_1+k_2\}$ and the probability of repair method $r$ being chosen is $\frac{w_{r}}{\sum_{z=k_1+1}^{k_1+k_2}w_z}$.

\subsubsection{Adaptive Weight Adjustment of Operator}

    We further explain how the weights $w_j$ can be automatically adjusted using statistics from earlier iterations.The basic idea is to keep track of the score for each heuristic, which measures how well the heuristic has performed after certain number of iterations. After the iteration $t$, we calculate the weight for heuristic $i$ to be used in iteration $t+1$ as follows:     
    \begin{align}
        w_{i,t+1} = \zeta w_{i,t} + (1-\zeta) \sigma, \sigma \in \{\sigma_1, \sigma_2, \sigma_3, \sigma_4\},\notag \\
        \zeta \in [0, 1], i \in \{1,2,...,k_2\}
    \end{align}
    where $\sigma_1$, $\sigma_2$, $\sigma_3$, $\sigma_4$ are score adjustment parameters when the founded candidate solution is a new best solution, a better solution,  a solution accepted or refused by acceptance criterion , respectively. And $\zeta$ is the reaction factor.

\subsubsection{Acceptance and Stopping Criteria}
    
    We incorporate three acceptance criteria: hill climbing acceptance criterion, record-to-record travel acceptance criterion \cite{dueck1993great}, and simulated annealing acceptance criterion \cite{kirkpatrick1983optimization}. Each of aforementioned criteria is accompanied by specific parameters to regulate its performance. We set the termination condition for AHGSLNS at a maximum of $iter^{max}$ iterations.

\subsection{Adaptive Survival Mechanism (ASM)}

   During the adaptive survival phase, we utilize PALNS to enhance the solution quality of individuals within the population in parallel and eliminate the least favorable ones, thereby adapting to different problem instances. Following the removal of inferior individuals from the population, the surviving individuals attain elite status, equipped with ALNS parameters and acceptance criteria well-suited for the specific characteristics of the current instance.
    
\subsection{Cooperative Evolution Mechanism (CEM)}

\subsubsection{Parent Selection and Crossover}

    The offspring generation scheme of AHGSLNS uses binary tournament to select two parents according to fitness, $P_1$ and $P_2$, and yield a single individual $C$. For MAVRP, we introduce both trip and journal crossover operators. The former exchanges trips within a journey, while the latter directly exchanges entire journeys. AHGSALNS incorporates the crossover operator when there is no improvement in the best solution of the population for $l_c * gen^{max}$, for $l_c \in [0,1]$ iterations, aiming to derive a promising solution structure.

\subsubsection{Diversification}

    For certain COPs encountered in real-world scenarios, mutating individuals becomes challenging due to the intricate nature of feasible solution structures. To overcome this limitation, AHGSALNS incorporates new individuals into the current population after $l_d * gen^{max} (l_d \in [0,1])$ iterations without enhancing the best solution within the population. This strategy aims to induce mutation within the population.

\section{COMPUTATIONAL RESULTS AND ANALYSIS}
\subsection{Computational Environment and Data Sets}
\paragraph{Computational Environment}

    The computational experiments related to the mathematical model and meta-heuristic algorithm were conducted on a computer equipped with an Intel Core $i5$-$7200U$ CPU and $8GB$ of RAM. The model was implemented in a general-purpose solver, GUROBI version $9.5.0$, with a maximum computational time of $3,600$ seconds.
    
\paragraph{Data Sets}

    The synthetic data sets are generated based on the random distribution (R) with the geographic type. The names of the synthetic instance graphs start with the geographic type, followed by the information on the number of depots, linehaul customers,backhaul customers, and vehicles. The dimension of the map sides is set to 30.

\subsection{Convergence Analysis}

    We configured the parameter for ALNS with  $iter^{max}$ set to 600, AHGSLNS with  $iter^{max}$  set to 20, $gen^{max}$ set to 30, $N$ set to 6, and $M$ set to 3. 
    This results in a total of 600 iterations for both ALNS and AHGSLNS. Our evaluation involves a comparison between mathematical formulations, ALNS and AHGSLNS algorithms on MAVRP instances (shown as Table 1). For each instance, the best known solution is given on column BKS and we define the computational gap as 
    \begin{align}
        GAP = \frac{MSP_{heuristic}-MSP_{formulation}}{MSP_{formulation}} \times 100\%
    \end{align}
    where $MSP_{algorithm}$ denotes the makespan of algorithm. We observed that both ILP and proposed algorithm were capable of generating the same BKS, which affirms the correctness and capability of all three methods in achieving optimal solutions on small instances. Furthermore, in comparison to ALNS, AHGSLNS exhibited superior average performance in terms of convergence on all instances of the MAVRP.

\subsection{Stability Analysis}

    To verify the stability of the algorithm, we introduce the coefficient of variation (CoV) and standard deviation (Std), which are important indicators to measuring the stability of a meta-heuristic. The definition of CoV is given as:
    % \begin{align}
    %     \mbox{CoV}=\frac{\frac{1}{T} \sum\limits_{i=1}^T \mbox{MSP}(x_i),}{\sqrt{\frac{1}{T}\sum\limits_{i=1}^T(\mbox{MSP}(x_i)-\mu)^2}}\times 100\%,
    % \end{align}
    \begin{align}
        \mbox{CoV}=\frac{\sigma}{\mu}\times 100\%,
    \end{align}
    \noindent where 
    \begin{align}
        \mu = \frac{1}{T} \sum\limits_{i=1}^T \mbox{MSP}(x_i),
        \sigma = \sqrt{\frac{1}{T}\sum\limits_{i=1}^T(\mbox{MSP}(x_i)-\mu)^2},
    \end{align}
    \noindent where $T=10$ refers to the number of runs, and $x_i$ refers to the makespan of the obtained best solution in the $i^{th}$ run. As shown in Table 1, the CoV value of AHGSLNS is lower than that of ALNS in up to $73.33\%$ instances which proves that the AHGSLNS is more stable than the ALNS method.

\subsection{Effect Evaluation of ASM and CEM}

    To evaluate the effect of the adaptive survival and cooperative evolution mechanism of AHGSLNS, we propose four ALNS with various parameter sets and acceptance criteria, namely A1, A2, A3, A4, B1 and B2 are AHGSLNS without cooperation evolution phase and AHGSLNS, respectively. See Table 2 for additional algorithm configurations while all other parameter settings remain consistent with Section 4.2.

    In Table 3, among  four different ALNS frameworks, the algorithms demonstrate varying degrees of success and shortcomings when compared to one another. Specifically, it is observed that A1 exhibits superior average performance in 60.00\% of the instances, whereas A2, A3 and A4 show 0.00\%, 26.66\% and 13.33\%, respectively. Additionally, for the best results of 10 runs, A1 shows better performance on 40.00\% of the instances, while A2, A3, A4 all display 26.66\%. The variation in performance is attributed to the fact that ALNS with different parameters and acceptance criteria may yield diverse results across instances. 

    Notably, B1 outperforms all for ALNS variants in up to 93.33\% of the instances, showcasing its efficacy in adapting to various scenarios of MAVRP. Further, compared with B1, B2 exhibits even greater advantages surpassing the results of B1 by 80.00\% in terms of average performance. This underscores the effectiveness of cooperative evolution mechanism. 
    
\section{Conclusion}
    In this paper, we study the MAVRP, a specific form of the COPs in the real world and propose a meta-heuristic algorithm named AHGSLNS, which adaptively combines the essential idea of ALNS and genetic algorithms, to address this NP-hard problem. Leveraging initial solutions generated by random constructive algorithm, we eliminate worse individuals considering various parameters and acceptance criteria. Additionally, through an adaptive survival mechanism, elite individuals were retained to adapt to the diverse attributes of MAVRP. Furthermore, by orchestrating crossover and diversification operators in an evolutionary manner, AGSALNS exhibited outstanding performance in computational adaptability, convergence, and stability across 30 test cases compared to traditional ALNS. This adaptive capability extends its applicability to solving other complex COPs problems, providing promising avenues for future research.

\appendix

% \section*{Ethical Statement}

% There are no ethical issues.

% \section*{Acknowledgments}

% Carton Transfer Unit robots warehouse system scenario (https://www.hikrobotics.com/en/mobilerobot/CTU) is supported by Hikvision in China. 

%% The file named.bst is a bibliography style file for BibTeX 0.99c
\bibliographystyle{named}
\bibliography{ijcai24}

\end{document}